\documentstyle[12pt,epsfig]{article}
\voffset - 1.0in
\title{INTRODUCTORY CALCULUS FROM THE VIEWPOINT OF 
NON-STANDARD ANALYSIS - DERIVATIVE OF SINE AND COSINE}
\author{JACK L. URETSKY \\ High Energy Physics Division, Argonne
National Laboratory\\(jlu@hep.anl.gov)}
\date{\today}
\begin{document}
\setcounter{page}{1}
\maketitle
\begin{abstract}
 This article exemplifies a novel approach to the teaching
of introductory differential calculus using the modern notion
of ``infinitesimal'' as opposed to the traditional approach using
the notion of ``limit''.  I illustrate the power of the new approach
with a discussion of the derivatives of the sine and cosine functions.
\end{abstract}
\section{OVERVIEW}
\indent Some teachers of mathematics may be surprised to learn 
that the concept of limits was a rather late (d'Alembert, 1754)
addition to the underpinnings of the calculus.  The inventors of
calculus worked directly with concepts such as ``infinitely little
lines'' (Newton, in Cohen and Westfall's {\em Newton}, Norton, New
York 1995, p.381), presumably adopting the most intuitively congenial
approach to the subject.  This fact suggests that it is
counter-intuitive to introduce students to the derivative as the
outcome of a limiting process. In the words of one mathematics teacher:
\begin{quotation} 
As someone who has spent a good 
portion of his adult life teaching 
calculus courses for a living, I can tell 
you how weary one gets of trying to 
explain the complex and fiddling 
theory of limits to wave after wave of 
uncomprehending freshmen. 
\end{quotation} 
(Rucker. {\em Infinity and the Mind}, Birkh\"{a}user, 
Boston 1982, p. 87).\\

\indent Newton and Leibniz, Leibniz especially, seem to have 
conceived of very small quantities, ``infinitesimals'', that 
lay outside of the ordinary number system.  Leibniz spoke 
of the derivative of a function as being the ratio of two 
such infinitesimals.  By the beginning of the twentieth 
century, however, the notion of ``infinitesimal quantity'' 
had disappeared from mathematics textbooks (although it 
persists in some physics and engineering texts as, for example,
in the ''principle of virtual work''.).  Then,
in 1966, Abraham Robinson showed that the real 
number system can be extended to incorporate the notion 
of infinitesimal in approximately the sense of Leibniz (Robinson,
 {\em Non-standard Analysis}, North-Holland, Amsterdam 1966).  
While the foundations and rules developed by Robinson 
involve advanced mathematics (a readable account may be found 
in Stroyan and Luxembourg, {\em Introduction to the Theory 
of Infinitesimals}, Academic Press, New York 1976) the results 
can be used to 
simplify greatly the solutions of many long-standing 
mathematical problems.

\indent There are existing texts that incorporate non-standard analysis
 \footnote{ {\em e.g.} H. Jerome Keisler, {\em Elementary Calculus, an 
Approach Using Infinitesimals}, Bodgen \& Quigley, 1971, (2nd ed. 
entitled {\em Elementary Calculus: an Infinitesimal Approach} 1986), 
and Keith Stroyan's {\em Calculus Using Mathematica} (1993) (second 
edition entitled {\em Calculus, the Language of Change} 1997)}  These 
have failed, I believe, to clarify the subject for undergraduates 
because the authors 
apparently felt honor-bound to clutter their presentation of the 
intuitively accessible ideas by first presenting rigorous 
justifications of the mathematics.  

\indent I contend that the basic ideas are easily understood 
by students who are probably totally impervious to the justifications. 
I also contend that students will experience far less discomfort
in dealing with infinitesimals from the outset than they 
presently experience while trying to absorb and utilize the notion 
of limits and related mathematical underpinnings.

\indent I show in this article how to avoid 
the early introduction of concepts that, while important to 
the rigorous development of the subject of 
calculus, are totally incomprehensible to beginning 
students.  The present approach is, instead, intended to 
build on students' intuitive concepts, with 
mathematical rigor generally left for advanced courses where 
rigor aids, instead of detracts from, the conceptual 
development. 

\indent  The underlying idea (for example in treating 
difference quotients in finding the slope of a curve), 
is to treat an extremely small increment of the 
independent variable (dubbed the `dibbl') as if its square were 
equal to zero.  This idea gives one an occasion to exercise 
students on the powers of large and small numbers, leading to a 
single new concept for the student to confront.  
This approach, which allows one to defer the notion of limits, 
can be justified using non-standard analysis.  The results can, 
of course, be duplicated using limits - often with much more work.

\indent The approach illustrated in this article differs in another 
significant way from that found in existing calculus textbooks.  
Almost all such textbooks seem to assume that a 
student learns the subject in a logically linear fashion. 
First we learn about functions and continuity, then about 
limits, then the derivative, and so forth.  The student 
isn't told, until somewhere around chapter 4 of the typical
text, that all of 
this mathematical apparatus is somehow related to solving 
problems concerning motion, velocity, acceleration, and 
related concepts. 
 
\indent  The present approach is structured 
to take into account that the learning of 
mathematics is a highly non-linear process.  Oftentimes 
the mathematical (or ``logical'') justification for a 
procedure cannot be understood by a student 
until long after the student has learned the procedure.  
All of us, I think, gain new insights into most topics 
in mathematics each time we revisit the topic.  
There is need, therefore, for texts that emphasize such 
revisitations as a part of the learning process.  

\indent Many introductory students are weak on 
algebra skills and have had minimal exposure to the 
concept of ``proof''.  This article incorporates exercises, 
and a pedagogical approach, intended to assist a 
teacher in addressing such weaknesses.  I try at all times 
 to incorporate the 
pedagogical principles set out in Arnold Arons' book {\em 
A Guide to Introductory Physics Teaching } (Wiley, New 
York 1990).  Although Arons' book ostensibly addresses 
the teaching of physics, much of the subject matter 
addresses means for remedying the mathematical 
deficiencies of typical engineering and science students. 
Also, I ordinarily prohibit, in keeping with this spirit, 
the use of calculators in working most (but not all) problems, 
accepting rough estimates of numerical values, where appropriate.

\indent I adhere, for the most part, to the
 American Standard Abbreviations for Scientific and 
Engineering Terms as reproduced, for example, in 
{\em Handbook of Chemistry and Physics\,} (63d Ed., CRC 
Press, Boca Raton, FL., 1982) 

\indent Section 2 of the present article is a highly abbreviated
introduction to differentiation, using the concept of the ''dibbl'',
a quantity so small that its square is zero.
The section omits, for the sake of brevity, many
topics that I would include in the first
six chapters of an introductory text such as the idea that
``velocity'' is the slope of the tangent to a
distance-time curve, or that ``rate'' is the slope of the tangent to
some other kind of curve. (See, {\em e.g.} excerpts from a proposed
text at http://www.hep.anl.gov/jlu/index.html.) It also
omits such vital topics as the product and chain rules.

\indent Section 3 uses the dibbl method to obtain the derivatives
of the sine and cosine functions.  I believe that the derivation
presented here is novel (for a textbook), and is much easier to
understand and reproduce than the presentations in other texts.
I deliberately omit, in that Section, any discussion of uniqueness.
``Less is more'', especially in introductory textbooks.

\indent  I simplify the language by speaking of the {\em slope of a
curve} when I mean the {\em slope of the tangent to the curve}.

\pagebreak

\section{ FINDING SLOPES}
\begin{quote} Your Majesty says, ``Kill a gentleman, and a
gentleman is told off to be killed.  Consequently, that gentleman
is as good as dead - practically, he is dead - and if he is dead,
why not say so? ''{\em The Mikado, Act II}, by Gilbert and Sullivan
\end{quote}

\indent I motivate the derivative, as pointed out in the introduction,
as measuring the rate of change of a quantity.  We deal with quantities
that can be represented by curves on a graph.  Our problem is to
find the rate of change of the quantity at a point on the graph.
That rate of change is measured by the slope of the tangent at the point
to the curve representing the function.

\indent In what follows, I shall introduce a quantity, that I
call  a ``dibbl'', so small that its square ``is as good as'' 
zero - practically, it is zero - and
if it is zero, why not say so?  The dibbl will make it easier 
for us to calculate slopes.

\subsection{NEW LANGUAGE}
\indent I show, in this Section, how to find the slopes of functions
that can be described by integer powers of an independent variable. The
problem is to find the equation of a straight line that is tangent to
the function at some point.  Once we know that equation we can just
read off the slope of the line.

\indent Let's try to find the slope of the tangent at some point, call it
$x_{1}$, of the quadratic equation $y(x) = cx^{2}$ where $c$ is some
constant.  Then the slope of a line joining that point to some other
point on the curve, call it $x_{2}$, is just the rise over the run so
that

\begin{eqnarray}
                \mbox{ slope} &= &\frac
{y(x_{2})-y(x_{1})}{x_{2}-x_{1}} \label{eq:q1}
 = \frac{cx_{2}^{2}-cx_{1}^{2}}{x_2-x_1} \nonumber \\ &=
&c(x_{2}+x_{1})
\end{eqnarray}

\indent We wanted to find the slope of the quadratic at a
point whose $x-$ coordinate is equal to
$x_{1}$.  But the process of taking the slope
involves two points because the slope is the rise of $y$ between two points
divided by the run along $x$ between the same two points.  So we have a
result that depends, not only on the point with $x=x_{1}$, where we
want the slope, but also on some other point with $x=x_{2}$.  We have,
in other words, the slope of a line that cuts the curve at two points.

\indent But we only want the tangent line to touch the curve at the one
point, namely where $x=x_1$. What to do?  Choose $x_{2}$=$x_{1}$!
Then Eq.~\ref{eq:q1} becomes:
\begin{equation}
\hbox{\rm slope} =2cx_{1}           \label{eq:q2}
\end{equation}

\indent Eq.\ref{eq:q2} shows that it is easy to find the slope of a quadratic function
$y(x) = x^{2}$ at a point $x = x_{1}$.  The slope is just the ratio of
the change in $y$ to the change in $x$ between two points $x_{2}$ 
and $x_{1}$ and then letting $x_{2} = x_{1}$.  
Writing the last statement in symbols,
we have
\begin{equation}
\mbox{slope} =  \frac{y(x_{2})-y(x_{1})}{ x_{2}- x_{1}}
 \label{eq:slope}
\end{equation}
and then let  $x_{2}= x_{1}$.
The same procedure should work for any smooth function $y(x)$, 
but the procedure
would become pretty tedious if the function is very complicated, as the
next exercises
demonstrate.

\begin{sf}
\begin{center}
{\bf EXERCISE}
\end{center}
\noindent{\bf 2.1}
Find the slope of the function $y(x) = x^{4}$ at the point
$x=x_{1}$
using the method of Eq. \ref{eq:slope}.  Show all your work.
Hint: (Check your result by noting that when $x_{1} = 3$ the value 
of the slope is 108.)

\noindent{\bf2.2}  You are running around a circular quarter-mile track.  
The coach tells you that your distance from the starting line is
increasing
proportionally to the cube of the time.  At the end of 1 second you have
traveled 3 m.\\
(a)  Write an equation for $s(t)$, the distance traveled, as a function of
$t$.  Check
your equation by noting that you must have traveled  24 m at the end of 2
seconds.\\
(b)  Your speed at any time $t$ is the rate of change of your distance at
time $t$.\\
Find your speed at any time $t$.  (Here we are writing just $t$, instead of
$t_{1}$).\\
(c).  A horse can gallop with a top speed of about 7 m/sec.  After how many
seconds
are you running at a horse's top speed?\\
(d)  Did you have to consider the shape of the track when you answered (a)
through (c)?\\
\end{sf}

\indent  It turns out that we can simplify the calculation of slope
 by taking the points $x_{2}$ and $x_{1}$ to be very close together.  This
makes sense, because at the end of the calculation we are going to let
$x_{2} = x_{1}$.  So let's start by simplifying the notation, dropping the
subscripts and restating the problem as follows:  {\bf Find the slope of the
function $y(x)$ at some value of $x$.}

\indent We will do this by
finding the slope of a straight line that intersects $y(x)$ at two points
that are
so close together that they are {\bf practically, but not quite, the same
point}.  We
write the two points as $y(x)$ and $y(x+dx)$.  Since $x$ and $x+dx$ are
nearly the same, we say that $x+dx$ is just a {\em{dibbl}} away from $x$.

\indent Also, since we know that $y(x)$ is described by a smooth curve, we
also know that the two points on the curve with y-values $y(x)$ and $y(x+dx))$
are just a dibbl away from each other.

\indent  I have already said that a dibbl is a quantity so small that
we can take its square to be zero ``for all practical purposes.''  So the
rest of this article makes use of ``the dibbl equation''
\begin{equation}
dx \times dx = dx^{2} = 0
\label{eq:dibbl}
\end{equation}

\indent  The dibbl equation comes into play when we recalculate
the slope of a quadratic at the point $y(x)$.
We write, for $y(x)=cx^{2}$:
\begin{eqnarray}
\mbox{slope} &= &\frac{y(x+dx)-y(x)}{(x+dx)-x} \label{eq:q3}\\
 & = &\frac{c(x^{2}+2cxdx +  \overbrace{dx^{2}}^{=0})-cx^{2}}{dx} \nonumber \\
 &= &\frac{2cx\times dx \hspace{-2ex}//}{dx \hspace{-2ex}//} = 2cx  \nonumber
\end{eqnarray}

Eq. \ref{eq:q3} no longer requires the step of setting the point $x_{2}$ equal
to the point $x_{1}$.  The dibbl equation makes that step unnecessary.

\indent  Note that the slope in Eq. \ref{eq:q3} is a function of $x$.  It
is, in fact,
a function of $x$ that was {\em{derived}} from the original function $y =
cx^{2}$. Mathematicians, for this reason, call the slope of a function y(x)
the {\em derivative}  and write it as a ratio of two dibbls, dy and dx.  To
summarize, then, for any smooth function $y(x)$ there is a slope function,
called the derivative,
\begin{equation}
\frac{dy}{dx} = \frac{y(x+dx)-y(x)}{(x+dx)-x} = \frac{y(x+dx)-y(x)}{dx}
\label{eq:deriv}
\end{equation}

\noindent{\bf IMPORTANT NOTE} Many textbooks, for historical reason, write the
derivative of a function $y(x)$ with a ``prime'' symbol.  So the symbols $y',
y(x)', y'(x), \mbox{and}\,\frac{dy}{dx}$ all mean exactly the same thing.  I
shall sometimes use the prime symbol in this article, so that students can
accustom themselves to both usages.
\begin{sf}
\begin{center}
{\bf EXERCISE}
\end{center}
\noindent{\bf 2.3}
Redo Exercise 2.1, this time using the dibbl method.
\end{sf}
\medskip

\subsection{POWERS OF SMALL NUMBERS;\\ THE DIBBL dx}
\indent We have introduced the dibbl as a quantity "that is so small that
its square is as good as zero".  "How small is that?", you might say.  So let's
get some experience with the squares of small numbers.

\begin{sf}
\begin{center}
{\bf EXERCISE}
\end{center}
\medskip\noindent{\bf 2.4}
You may use your calculator for this exercise: \\
(a)  How much is the square of $\frac{1}{5}$?  (Hint:  If you didn't get
$\frac{1}{25}$, ask your teacher to send you back to third grade).\\
(b)  Now take a piece of paper and a ruler marked in centimeters and
draw a straight line that is 10 cm long.  Divide the line into 5 equal parts
by marking little perpendicular lines ("tick marks") on the line, like the
longer vertical lines in Fig. 1.\\
\indent Next, divide each division into 5 equal parts, again using tick marks.
Into how many divisions have you divided the 10 cm line?  Your line should
resemble the line in Fig. 1.\\
(c)  Does the distance between the little tick marks in Fig 1 represent the
square of the distance between the large tick marks?  (Hint:  If you answered
"No", go back and redo part (a) of this Exercise).  Now complete the following
table, using your calculator if you so desire:\\
\begin{figure}
\centerline{\epsfig{figure=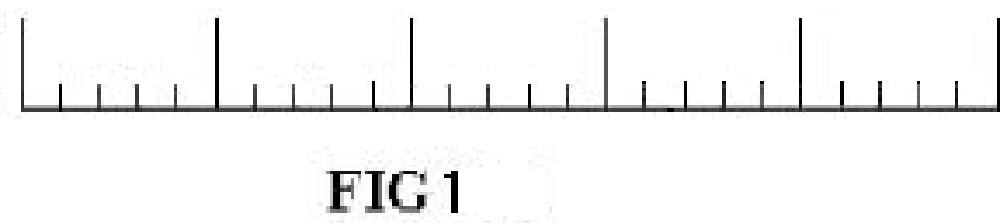, height=1.11in}}
\end{figure}
$\begin{array}{clcrcl}
  x= &  \frac{1}{5} =&  .2  \qquad &  x^{2}  = & \frac{1}{25}  &  = .04  \\
   = &  \frac{1}{50} =&     &            &               & =     \\
   = &  \frac{1}{500} =&   &   &     & =     \\
   = &  \frac{1}{5000} =&   &   &     & =     \\
   = &  \frac{1}{50000} =&   &   &     & =
\end{array}$\\
(d)  Are you tired of writing out the zeros in part (c)?  Let's do the same
calculations in power-of-ten notation, which we will use in all future
exercises that involve numbers:\\
$\begin{array}{clcrcl}
 x=&  \frac{1}{5} =&2 \times 10^{-1}  &x^{2}  = & \frac{1}{25} &= 4 \times 
10^{-2}   \\
  &  \frac{1}{5 \times 10} = &   &   &     & =    \\
  &  \frac{1}{5 \times 10^{2}} = &   &   &     & =     \\
  &  \frac{1}{5 \times 10^{3}} = &   &   &     & =   \\
  &  \frac{1}{5 \times 10^{5}} = &   &   &     & =
\end{array}$\\
\end{sf}

\indent  Exercise 2.3 shows that the square of a fraction means ``taking a
fraction of a fraction".  If the fraction is a very small number then its
square may be very small indeed.  How small must the number be so that the
square ``is as good as" zero?  In other words, how small is a dibbl?  This 
is a profound question.  It suffices to say for now that a dibbl is not an ordinary number, but is something that is smaller than any fraction!

\indent We close this subsection by reminding ourselves why it is
difficult to assign a value to $\frac{1}{0}$.

\begin{sf}
\begin{center}
{\bf EXERCISE}
\end{center}
\noindent{\bf 2.5}\\
(a) Ask your calculator to find a value for $\frac{1}{0}$, and write
down the result.  If you got a value, you are finished.\\
(b)  If your calculator couldn't do the problem try dividing $1$ by
a lot of values in the denominator that are close to zero.  Then you
can try to guess the result.  One way to do this is to take a
sheet of graph paper, and plot values of $\frac{1}{x}$ for
a lot of different fractional values of $x$.  (Don't forget to
include negative values).  The answer should be roughly
half way between the value for the smallest positive value of
x and the smallest negative value.  What is your result?\\
(c)  If you are in a study group, compare your result 
with others in your group.  Can
you all agree on the same value?\\
(d)  What would you say to someone who claims that 
division by zero is ``undefined''?  
Explain your answer, referring to your results in part (b).\\
\end{sf}
\medskip
\subsection{$dx^{2}$ IS ZERO SIMPLIFIES\\ THE BINOMIAL THEOREM}
\begin{sf}
\begin{center}
{\bf EXERCISE}
\end{center}
\end{sf}
\noindent{\bf 2.6}
\begin{sf}
(a).  Write down the derivative functions (the slope function, remember?)
$\frac{dy}{dx}$ for $y(x) = x, x^{2}, x^{3}, x^{4}$.  (You worked these
out in Subsection 2.1).  \\
(b)  Guess what the derivative function is for $y(x) = x^{n}$ where $n$ stands
for any positive integer (do not use a specific numerical value for  $n$), and
write down your answer.
\end{sf}
\medskip

\indent We now help you work out the answer to Exercise 2.5(b), using the
 definition of the derivative in Eq. \ref{eq:deriv}.  The first step is to
expand
the product $(x + dx)^{n}$, remembering that all powers of $dx$ higher than
the first power are zero.  We write:
\begin{equation}
y(x+dx) = (x + dx)^{n} =\overbrace{(x + dx)(x + dx) \cdots
(x + dx)}^{\mbox{n\, times}} \label{eq:binom}
\end{equation}
\indent The highest power of $x$ in Eq. \ref{eq:binom} comes from
multiplying all of the $x's$ together.  There is only one way to do this,
 so there is a leading term $x^{n}$.

\indent The next term will come from multiplying together 
$(n-1)$ $x's$  and one $dx$.  How many such terms are there?  
That's an easy question, because each term omits one $dx$ and 
there are exactly $n$ different $dx's$ in the product in
Eq. \ref{eq:binom}.  The second term is therefore $nx^{n-1}dx$.

\indent What about the next term?  It will have the product of
$(n-2)$ $x's$, two $dx's$ and a numerical coefficient, let's call it $a$.
  There will be more terms, involving products of three and even more
 $dx's$ (if $n$ is bigger than 2).  The numerical coefficient that goes
with each term is known from the {\em binomial theorem} that people
 learn about in pre-calculus courses.

\indent But we don't care about the value of $a$ and all those
numerical coefficients because they all involve products of two or more
$dx's$, and since $dx$ is a dibbl the product
of two or more $dx's$ is zero.  The result, then, is
$(x + dx)^{n} = x^{n} + nx^{n-1}dx$, which is the statement of the binomial
theorem  for $x$ plus a dibbl taken to any power $n$.\\

\indent Eq \ref{eq:deriv} gives, for $y(x) = x^{n}$:
\begin{eqnarray}
\frac{dy}{dx}& = &\frac{d(x^{n})}{dx} \nonumber \\
             & = &\frac{(x^{n} + nx^{n-1}dx) - x^{n}}{dx} \nonumber \\
             & = &\frac{nx^{n-1}dx}{dx} = nx^{n-1} \label{eq:power}
\end{eqnarray}
Eq \ref{eq:power} is often called "the power law".

\begin{sf}
\begin{center}
{\bf EXERCISE}
\end{center}
\noindent{\bf 2.7}
Calculate the derivatives $y' \equiv \frac{dy}{dx}$ ($\equiv$ means that
the two expressions mean the same thing), using the method of
Eq \ref{eq:power} when\\
(a) $y(x) = x^{9}$; \\
(b) $y(x) = 5x^{17}$; \\
(c) $y(x) = 6x^{5} + 4x^{4}$\\
Show all of your work.
\end{sf}
\medskip
\subsection{``DERIVATIVE'' IS A FANCY NAME\\ FOR VELOCITY OR RATE OF CHANGE}
\begin{sf}
\begin{center}
{\bf EXAMPLE 2.1}
\end{center}
Find an expression for the straight line that is tangent to the curve described
by $y(x) = \frac{1}{7}x^{5}$ at the the point $x = 2$.  Show your work and
graph the
result. \\
\noindent{\bf Solution:}\\
(1)  Let $y_{t}(x) = a + bx$ be the equation for the straight line.  Then
our job is to determine the values of the constants $a$ and $b$. \\
(2)  Since $y_{t}(x)$ is a tangent, it must have the same slope as $y(x)$
at $x = 2$.
But, as you learned in algebra, the slope of $y_{t}(x)$ is $b$.  
The slope of $y(x)$ at any point $x$ is $\frac{dy(x)}{dx}$\\
(3) At $x=2$, $\frac{dy(x)}{dx}= \frac{5}{7}x^{4} =
\frac{80}{7}$.  So $b = \frac{80}{7}$. \newline
(4)  Also, $y_{t}(x)$ and $y(x)$ must be equal at the point of tangency, since
they share the same point.  So $ a + \frac{80}{7}x = \frac{1}{7}x^5$  at $x
= 2$, or
$a + \frac{160}{7} = \frac{32}{7}$.  Thus, $a = - \frac{128}{7}$.  \newline
(5) The solution is $y_{t}(x) = \frac{-128 + 80x}{7}$.  Fig. 2 is a graph of
the solution showing the tangent line and the curve.
\end{sf}
\medskip
\begin{figure}
\centerline{\epsfig{figure=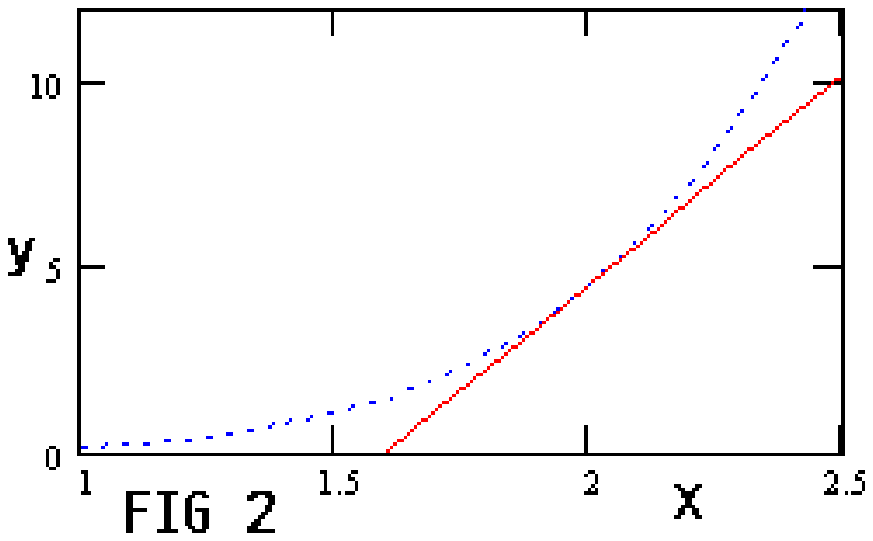, height=2.0in}}
\end{figure}

\indent Recall why we are interested in the slopes of tangents to curves.
We began our discussion of calculus by demonstrating that
``the slope of a line on a
distance-time graph \ldots corresponds to velocity".  So if we are
given a function that relates distance to time, the derivative
of that function, being the the slope of the function, corresponds
to the velocity as a function of time.  The next exercise makes use
of that relationship.

\begin{sf}
\begin{center}
{\bf EXERCISE}
\end{center}
\noindent{\bf 2.8}
(a)  A rock is hurled upward from a platform that is $H$ meters ($Hm$) above
the ground with an initial velocity $v_{0}\mbox{ \(m/sec\)}$.  We learn in
physics that the height of the rock above the ground at time $t$ after it
left the platform is
\begin{equation}
h(t) = H + v_{0}t - ct^{2} \label{eq:vert}
\end{equation}
where $c$ is a constant.
Find the velocity function $v(t)$ of the rock.  (Hint: use the
method of the previous subsection.  Make sure that you treat "h" and "H"
as two different symbols.)\\
(b)  Find an expression for $t_{m}$, the time when the velocity
function equals zero.\\
(c) Find an expression for the slope of $h(t)$ at time $t_{m}$.\\
(d) Find an expression for the straight line that is tangent at
time $t_{m}$ to the curve described by $h(t)$ on a distance-time graph.\\
(e)  Find an expression for the maximum value of $h(t)$ as a function of $t$.
  (Hint: to check your work, the answer to this part, when
$H = 8\mbox{ \(m\); }$ $v_{0}= 12\mbox{ \(m/sec\);}$ and
$c = 5\mbox{ \(m/sec^{2}\)}$, is $15.2\mbox{ \(m\)}$.\\
(f) Make a graph of $h(t)$ on the vertical axis against $t$ on 
the horizontal axis, using the values in part (e).  
The horizontal axis ofyour graph should be at least
15 cm long, and the horizontal scale should range 
from 0 to 4 secs.  Calculate $h(t)$ for at least 8 values of $t$.\\
\end{sf}

\noindent{\bf 2.9}
\begin{sf}
Now let's work out the equation for the tangent line to any power function.  So let
$y(x)=cx^{n}$ where c can be any constant and n is any positive
integer.  Write the tangent line, as before, $y_{t}(x) = a + bx$.\\
(a) Find the value of $b$ at an arbitrary point $x_{0}$.  (Hint: b
will be a function of $c, n, and x_{0}$).\\
(b) Now find the value of $a$ and write your expression for $y_{t}(x)$.\\
(c) Check your result by letting $c$ and $n$ take the values
used in the example.\\
\end{sf}

\noindent{\bf 2.10}
\begin{sf}
And, for an algebra brush-up, solve for t by completing the square
$ at^{2} + bt +c = 0$.  Show all of your work.
\end{sf}

\pagebreak
{\bf I OMIT, FOR THE SAKE OF BREVITY, SECTIONS ON THE PRODUCT AND 
CHAIN RULES}

\section{TRIGONOMETRIC FUNCTIONS}
\begin{quote} ...they be crammin my mind with [stuff] 
they want us to remember, ... all kinds of crap like 
coaxiel coordinates, cosine computations, spheriod 
trigonometry, ...
\end{quote}
\begin{flushright}
-Winston Groom ({\em Forrest Gump})
\end{flushright}

\indent In this section I present a novel (for an
introductory textbook) discussion of the derivatives
of the sine and cosine functions.  Many students at
this level have not learned to view the trigonometric
functions as functions. I therefore begin by introducing
that viewpoint.  "Spheriod trigonometry", however, can
safely be left for a later course.

\subsection{Circular Motion}
\indent The sine and cosine can be used in the description
\begin{figure}
\centering \includegraphics{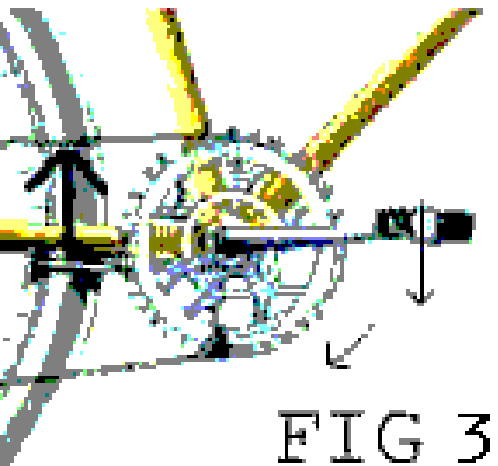}
\end{figure}
of circular motion.  Look, for example, at Fig 3
which is a close-up view of the working part of a bicycle.
The bicycle moves forward as the rear wheel is forced to 
rotate by the chain, which is driven by the rotation of
the sprocket.  But the sprocket is made to rotate by 
the up-and-down motion of the bicyclist's feet.

\begin{figure}
\centering \includegraphics{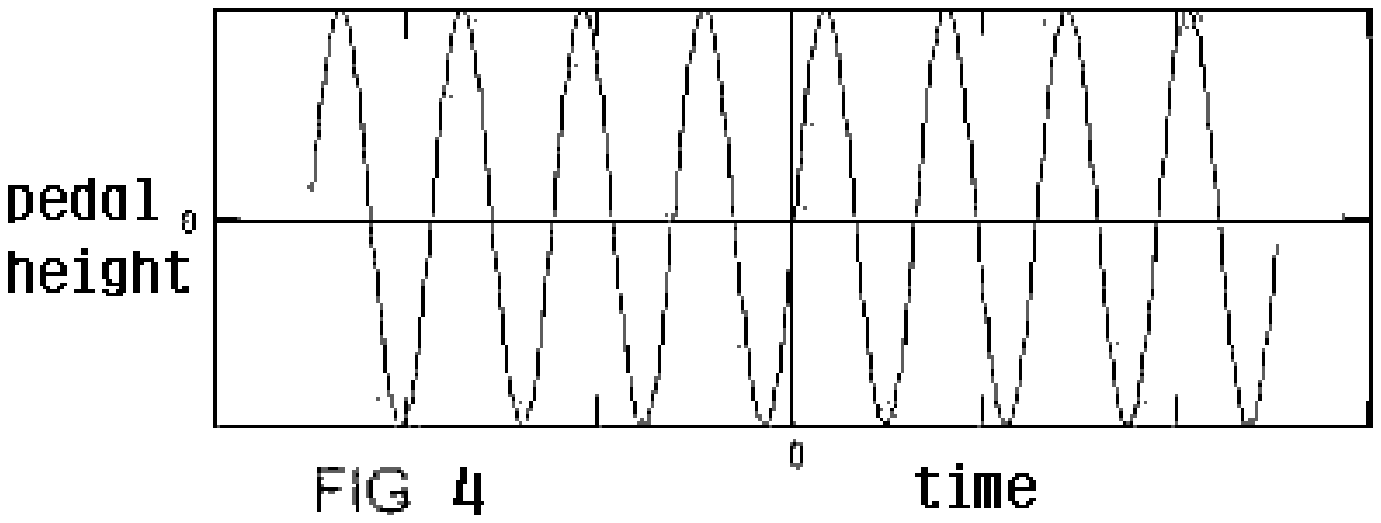}
\end{figure}
\indent Fig 4 is a graph that plots the vertical 
position of the bicyclist's left foot against time, 
assuming that the bicyclist rotates the sprocket 
with constant speed.  The vertical scale shows the 
height of the bicyclist's left foot above or below 
the position shown in Fig. 3, which is taken to be 
at zero time.

\indent The curve shown in Fig. 4 represents the sine
function that you learned about in connection with right 
triangles.  Its connection with circular 
motion should become clear after the next exercise.
\pagebreak
\begin{sf}
\begin{center}
{\bf EXERCISE}
\end{center}
\nopagebreak
\noindent{\bf 3.1}\\
(a) The circle shown in Fig. 5 has a radius that is exactly 
equal to one unit of length.  Express the distances $X$ and $Y$, 
shown in the figure, in terms of the sine and cosine of the
angle $\theta$.  (Hint: if you have trouble with this part,
review the sine and cosine definitions in one of your old math books - 
or a dictionary.)\\
(b) Make a graph of the length $Y$ as a function of the 
angle $\theta$, using $\theta$ intervals of $10^{\circ}$ 
between $0$ and $360^{\circ}$.  Does your graph resemble 
a portion of Fig. 4?  Be sure to plot the $Y$ values 
on the vertical axis.  \\
\begin{figure}
\centering \includegraphics{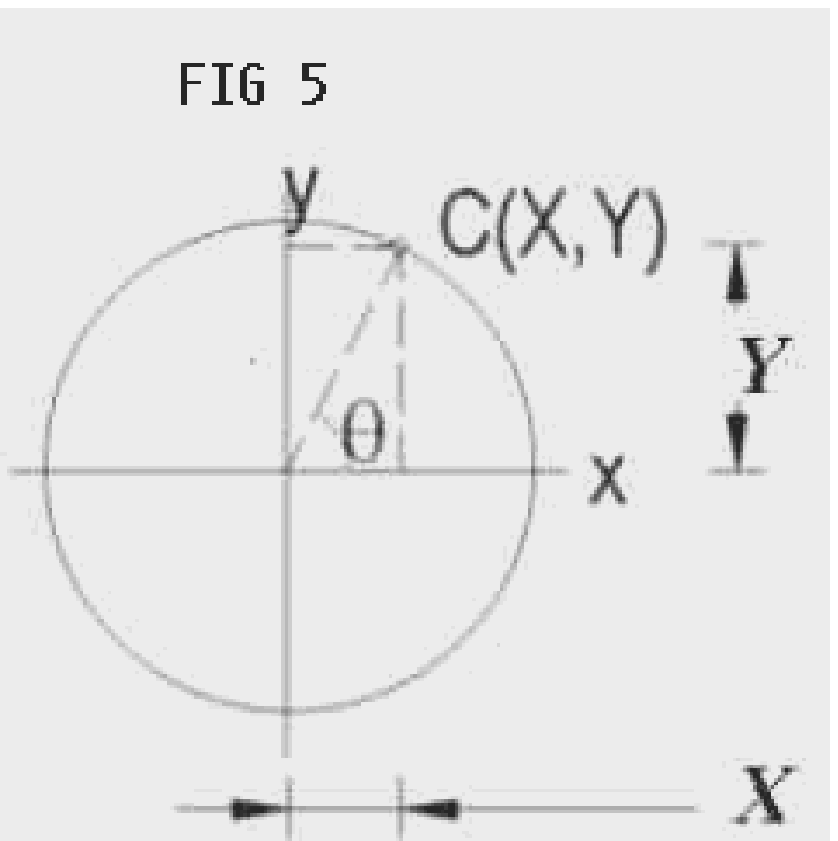}
\end{figure}
(c) Suppose the point $C$ in the figure moves around the 
circumference of the circle with constant speed so that 
$\theta = pt$, $t$ being the time.  Can you now make a 
graph of the length $Y$ as a function of the time?  
(Hint: Can you save yourself a lot of work be just 
relabeling the x-axis of your previous figure?  Does
it help to use time units of $\frac{degrees}{p}$, even
though you don't have a numerical value for $p$?\\
(d) Now sketch, without calculating, the graph of 
$Y$ for one more complete revolution of the point 
$C$ around the wheel.  Compare your graph with Fig. 3.  
Explain in one or two sentences what your graph 
would look like if it showed $Y$-values for many
revolutions of the wheel.\\
\end{sf}

\indent The last exercise demonstrates that the 
vertical position of the point 
$C$ in Fig. 5 can be written as 
a function of time (assuming that the sprocket 
rotates with constant speed) as 
\begin{equation}
Y(t) =R\sin (pt)
\label{eq:Y(t)7}
\end{equation}
where $p$ is the rotational speed 
of the sprocket and $R$ is the radius of the circle
described by the pedal.  Similarly the horizontal 
position of $C$ would be written as 
\begin{equation}
Y(t) =R\cos (pt)
\label{eq:X(t)7}
\end{equation}
The vertical positions of
the pedals in Fig. 3 would be equal to 
$Y(t)$, as Fig. 4 indicates.\\
\indent Let's get more familiar with the sine function.

\begin{sf}
\begin{center}
{\bf EXERCISE}
\end{center}
\nopagebreak
\noindent{\bf 3.2}\\
Suppose that a function $y(t)=R \sin (pt)$, where 
R and p are constants.  You may use Fig. 4 as a 
graph of $y(t)$.\\
(a) Suppose that $pt$ is equal to some fixed number 
$\theta$ that lies between $0^{\circ}$ and $360^{\circ}$.
How often will $y(t)$ repeat the same value of $y$ as 
$t$ increases forever?\\
(b) How many degrees separate the points on a graph of 
$y(t)$ as a function of $pt$ where y(t) not only has the
same value, but has positive slope?  (Hint: How many 
degrees of rotation bring the point $C$ back to the 
position shown in Fig. 5?)\\
(c) What are the values of $y(t)$ at the points where 
$y(t)$ has zero slope?\\
(d) Suppose that $x(t)=R \cos (pt)$.  What are the values
of $x(t)$ at the points where $y(t)$ has zero slope?\\
(e) What are the values of $y(t)$ at the points where 
$x(t)$ has zero slope?\\
(f) From the fact that $\cos (\theta) = \sin 
(\theta + 90^{\circ})$, what can you say about the 
{\em shape} of the graph of $x(t)$?  (Hint: could you 
match the graphs of $x(t)$ and $y(t)$ by 
sliding one over the other?)\\
(g) Sketch the derivative of the curve shown in 
Fig. 3.  Where the curve is increasing, the 
derivative is positive, the faster the increase 
then the more positive the derivative; 
correspondingly for places where the curve is 
decreasing.  Compare your curve with the cosine 
curve of part (f).  Redo your sketch after you 
have finished this subsection.
\end{sf}

\indent The functions $x(t)$ and $y(t)$ in the last 
exercise keep repeating themselves as $t$ increases.  
The repetitious behavior is not surprising, because it 
is connected with the motion of a point that is moving 
around the same circle, over and over again, with 
constant speed.

\indent We close this subsection with two little exercises 
on trigonometry.

\begin{sf}
\begin{center}
{\bf EXERCISES}
\end{center}

\noindent{\bf 3.3}\\
(a) A flagpole that has a circumference 
of $w$ inches is planted on level ground.  You are standing 
$H$ feet away from the center of the base of the flagpole.  
You find that the line-of-sight from yourself to the top of 
the flagpole makes an angle of $\theta$ degrees with the 
ground.  Your eyeballs are $h$ feet above the ground.  Find 
an expression for the height of the flagpole above the ground.\\
(b) Check your answer by noting that when $\theta =30^{\circ}$,
$H =200$ and $h =5$, the height is about $120$ feet.\\
(c) How much does your answer change if the circumference 
is doubled?\\
\noindent{\bf 3.4}\\
Recall that the addition formula for the sine function is 
\begin{equation}
\sin (\alpha + \beta )=\sin (\alpha ) \cos (\beta ) + 
\cos (\alpha ) \sin (\beta )
\label{eq:7add}
\end{equation}
Suppose that the symbol $\theta$ stands for some fixed angle.  
List all of the angles $\theta +  \alpha$ for which the graph of the 
$\sin$ function looks exactly the same as it does at $\theta$.
(Hint: a possible wrong answer would be
$\theta + n \times 180^{\circ}$, where $n$ is any 
positive or negative integer.)\\
\end{sf}
\subsection{DERIVATIVE OF THE PYTHAGOREAN THEOREM}
\subsubsection{Derivative of the Sine and Cosine}
\indent You showed, in part (a) of Exercise 3.1 that
the distances $X$ and $Y$ are respectively equal to 
$\sin \theta$ and $\cos\theta$.  The radius of the 
circle in Fig. 5, you were told, 
is 1 unit so, by the Pythagorean 
theorem, $X^{2} + Y^{2} =1$.  Expressing $X$ and $Y$ in 
terms of $\sin \theta$ and $\cos \theta$ gives the 
familiar relation from trigonometry \\
\begin{equation}
\sin ^{2}\theta + \cos ^{2}\theta = 1
\label{eq:s-sqd}
\end{equation}
\begin{sf}
\begin{center}
{\bf EXERCISES}
\end{center}
\nopagebreak

\noindent{\bf 3.4}\\
Suppose that the radius of the circle in Fig. 5 
is different from 1 - call it $r$, where $r$ is 
some fixed length.  What would then be the sum of 
the squares of the $\sin$ and $\cos$?  (Hint: divide 
the Pythagorean theorem by $r^{2}$.)\\

\noindent{\bf 3.5}\\
Take the derivative of both sides of the equation 
\begin{equation}
\sin ^{2}\theta + \cos ^{2}\theta = 1
\label{eq:sinsq7}
\end{equation}
Since you don't know the derivative of either $\sin$ 
or $\cos$ you may respectively write these as $\sin ' \theta$ 
and $\cos ' \theta$.  
(Hint: use the chain rule and the product rule).\\
\end{sf}

\indent The derivative of the Pythagorean Theorem (be sure
to do the last exercise),
expressed in terms of the $\sin$ and $\cos$, can be 
rearranged so that it reads
\begin{equation}
\frac{\sin ' \theta}{\cos \theta} = 
-\frac{\cos ' \theta}{\sin \theta}
\label{eq:difeq7}
\end{equation}
This equation contains two unknown functions, $\sin 
' \theta$ and $\cos ' \theta$.  Equations 
of this sort, containing unknown derivative functions 
are called {\em differential equations}.  They are 
often solved by guessing the solutions, as we now 
illustrate.

\indent Eq \ref{eq:difeq7} is an equality between two 
ratios, each involving an unknown derivative function.
One possible guess is that each ratio is equal to a 
constant, call it $A$.  If this were true then each derivative would
itself be proportional to a trigonometric function.  
Would this be a sensible result?

\indent We know that both the $\sin$ and $\cos$ are 
functions that keep repeating themselves at regular 
intervals (intervals of $360 ^{\circ}$; such functions 
are called {\em periodic} functions with period $360 ^{\circ}$).
  We also know 
that the {\em slope} of the sine function is zero at 
each point where the {\em value} of the cosine function 
is zero, and {\em vice versa}.

\indent It is obvious that the slope function for the 
$\sin$ is a periodic function that has a zero value at 
each point where the slope of the $\sin$ is zero.  Surely 
the $\cos$ is such a function.  So our guess is that 
\begin{equation}
\sin ' \theta \equiv \frac{d\sin \theta}{d\theta} =
A\cos \theta
\label{eq:dsin7}
\end{equation}
where $A$ is some constant.
\begin{sf}
\begin{center}
{\bf EXERCISE}
\end{center}

\noindent{\bf 3.6}\\
(a) Obtain an equation for $\cos ' \theta$.  (Hint:
you should be able to make use of the argument in the
preceding paragraphs.)\\
(b) Use the discussion of Exercise 3.2(f) to argue, 
in two or three sentences, that the guesses for 
$\sin '$ and $\cos '$ are reasonable.  Use 
sketches to support your argument.\\
\end{sf}

The guess works!  Solutions to the differential equation 
are given by Eq \ref{eq:dsin7} and 
\begin{equation}
\cos ' \theta \equiv \frac{d}{d\theta}\cos 
\theta = -A\sin \theta
\label{eq:dcos7}
\end{equation}

But what about that constant $A$?

\indent The first clue as to the meaning of the constant
is that we have never, in the preceding equation, specified 
the units that we are using to measure the angle $\theta$.
We will explore this point further in the next subsubsection, 
taking advantage of the fact that your calculator can 
calculate trigonometric functions in different kinds of units.
\subsubsection{The Constant A}
\begin{sf}
\begin{center}
\nopagebreak
{\bf EXERCISE}
\end{center}

\noindent{\bf 3.7}\\
This is a calculator exercise.
Estimate the slope of the sine function at $0^{\circ}$ using 
two different sets of units, and using Eq. \ref{eq:slope} of 
Section 2, as follows:\\
(a) Take $x_{2}= 5^{\circ}$, $x_{1}=0$.  Make certain that 
your calculator is set to calculate in degrees.  What is the 
value of $A$ for this case (please, no more than 3 significant 
digits)?\\
(b) Engineers sometimes measure angles in units called grads.  
A grad is $\frac{1}{400}$ of a circle, so that $5^{\circ}= 
5.55$ $grads$.  Set your calculator on ``grads'' and repeat 
the calculation of part (a) to obtain an estimate of the value 
of $A$ for the grad units.\\
(c) In part (b), did the value of the sine function change when 
you changed units?  Can the value of the sine of an angle 
depend on the units that you use to measure the angle?  
(Hint: recall the definition of the sine as the ratio of 
lengths of a side and the hypotenuse of a right triangle).
\end{sf}

\indent If you did your calculation correctly, you should have 
found that slopes of the $\sin$ and $\cos$ functions are about 
$10\%$ bigger (in magnitude) in degree units than 
they are in grad units.  So why don't we simplify matters 
by choosing an angular measure for which $A=1$?

\indent What kind of units for measuring angles give $A=1$?
We can investigate this question by looking at the sine of 
a very small angle, namely, an angle that is a dibbl.  
\begin{sf}
\begin{center}
{\bf EXERCISE}
\end{center}

\noindent{\bf 3.8}\\
(a) Write out the definition of the derivative 
in the expression from Eq \ref{eq:dsin7}, with $A=1$, 
\[\frac{d\sin \theta}{d\theta} 
\equiv \frac{\sin (\theta +d\theta)-\sin \theta}{d\theta}
 = \cos \theta,\]
but set $\theta$ equal to $0$ on both sides of the 
equation.\\
(b) Solve the equation that you obtained in part
(a) for $\sin (d\theta)$ to obtain an expression for
the sine of a very small angle.  Be sure to use the 
correct values of $\sin$ and $\cos$ at $\theta =0$.
Note that the $\sin$ of a dibbl is proportional to a
dibble, and, therefore, is also a dibbl.\\
(c) Now find $\cos d\theta$ ({\bf Hint}: recall that the 
square of a dibbl is zero).\\
\end{sf}

\indent The sine of a very small angle, expressed in
units where $A=1$
is just equal to the angle itself, as you have just 
demonstrated. That is, $\sin (d\theta) =d\theta$ in
this special system of angular units.
What are these special units?

\indent Since, as you have also shown, the cosine of 
a dibbl angle is just $1$, the side adjacent the angle
is equal to the hypotenuse.  The side opposite, 
therefore, begins and ends on the arc of a circle of 
radius equal to the hypotenuse.  

\begin{figure}
\centering \includegraphics{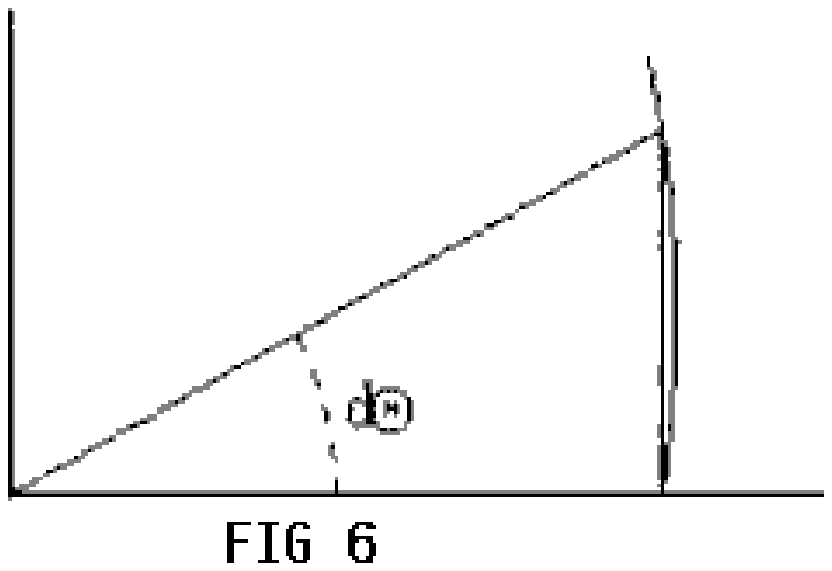}
\end{figure}
\indent In other words, the length of the side opposite the 
dibbl angle is just equal to the length of the arc 
swept out by the angle, as is shown in Fig. 6 (which 
exaggerates the curvature of the arc). 
But the $\sin$ is the length of the side opposite divided 
by the hypotenuse, which is the radius of the arc.

\indent Conclusion: the angular measure corresponding to 
$A=1$ is the arc length divided
by the radius.  You should recognize this as the definition
of radian measure.

\indent We end with an exercise that invites you to practice 
some of the skills you have developed thus far.
\begin{sf}
\begin{center}
{\bf EXERCISE}
\end{center}

\noindent{\bf 3.9}\\
Find the slope function for the function
\begin{equation}
F(t) = (\frac{t^{\frac{5}{3}}}{5 + 6t^{\frac{5}{3}}})^
{\frac{2}{7}}
\end{equation}
\end{sf}
Check your answer by finding the slope at $t=2$ where 
the slope is (to 3 significant digits) $.0134$.
\subsection{...AND, IN CLOSING}
\indent Trigonometry was important to early astronomers 
who, we blush to say, sometimes received income for giving 
astrological predictions.  The sine of an angle had probably 
received its present definition by about the fifth 
century A.D.  The Hindu mathematician \~{A}ryabhata 
(born in 476) is credited with this development.

\indent The study of the $\sin$ and $\cos$ as mathematical 
functions grew when the development of the calculus 
excited the interests of 17\raisebox{.7ex}{\em th} and 
18\raisebox{.7ex}{\em th} century mathematicians.  The 
Swiss mathematician Leonhard Euler, of whom you will hear 
much more in your studies of mathematics, may have been 
the first to treat the sine and cosine as mathematical 
functions and knew their derivative functions.

\indent Euler, in fact, wrote the derivative of the $\sin$ 
in the form of a differential equation
\begin{equation}
\frac{dy}{dx} = \frac{1}{\sqrt{1-x^{2}}}
\label{eq:arcsin7}
\end{equation}
You don't recognize it?  Invert each side, substitute an 
angle $\theta$ for the variable $y$ and $\sin \theta$ for 
the variable $x$.  A little algebra will give you the 
result of this subsection, for the case where $\theta$ is 
measured in radians.

\indent Now redo Exercise 3.2(g).
\section*{Acknowledgements}
I am indebted to Professor David Tartakoff for an explanation 
of non-standard analysis and for many helpful suggestions.  I 
am also indebted to Professor Gordon Ramsey and my colleagues
Brian Harris and Zachary Sullivan for their comments
on an early version of this manuscript.

This work was supported in part by the US Department 
of Energy, High Energy Physics Division, under Contract 
W-31-109-Eng-38. 
\end{document}